\documentclass{commat}

\usepackage[shortlabels]{enumitem}
\allowdisplaybreaks

\title{%
    Proper biharmonic maps on tangent bundle
}

\author{%
    Nour Elhouda Djaa, Fethi Latti and Abderrahim Zagane
}

\affiliation{
    \address{Djaa N.E.H --
    Department of Mathematics,
Relizane University, Algeria.
}
    \email{%
    Djaanour@gmail.com
}
    \address{Latti F. --
     Department of mathematics, SALHI Ahmed University Naama, Algeria.
}
    \email{%
    lattifethi2017@gmail.com
}
    \address{Zagane A.--
    Department of Mathematics,
Relizane University, Algeria.
}
    \email{%
    zagane2016@gmail.com
}
}

\abstract{%
     this paper, we define the Mus-Gradient metric on tangent bundle $TM$ by a deformation non-conform of Sasaki metric over an n-dimensional Riemannian manifold $(M, g)$. First we investigate the geometry of the Mus-Gradient metric and we characterize a new class of proper biharmonic maps. Examples of proper biharmonic maps are constructed when all of the factors are Euclidean spaces.
}



\VOLUME{31}
\YEAR{2023}
\NUMBER{1}
\firstpage{137}
\DOI{https://doi.org/10.46298/cm.10305}

\begin{paper}

\section{Introduction}\label{intro}

Let $\varphi:(M^{m},g)\rightarrow(N^{n},h)$ be a smooth map between two Riemannian manifolds. Such $\varphi$ is said to be harmonic if it is a critical point of the energy functional
\begin{equation}\label{Eq0.1}
    E(\varphi;D) = \frac{1}{2}\int_{D}|\tau(\varphi)|^{2}\,v_{g},
\end{equation}
for any compact domain $D\subseteq M$. Equivalently, $\varphi$ is harmonic if it
satisfies the associated Euler-Lagrange equations given as follows :
\begin{equation}\label{Eq0.2}
\tau(\varphi) = Tr_{g}\nabla d\varphi = 0.
\end{equation}
Here $\tau(\varphi)$ is the tension field of $\varphi$. We refer one to~\cite{MC.D},~\cite{J.L1},~\cite{J.L2},~\cite{KZD} for background
on harmonic maps. As a generalization of harmonic maps, biharmonic maps are
defined similarly, as follows:

 A map $\varphi$ is said to be biharmonic if it is a critical
point of the bi-energy functional
\begin{equation}\label{Eq1}
    E_{2}(\varphi;D) = \frac{1}{2}\int_{D}|\tau(\varphi)|^{2}\,v_{g},
\end{equation}
over any compact domain $D$. Equivalently, $\varphi$ is biharmonic if it satisfies the associated Euler-Lagrange equations
\begin{align}\label{Eqe1}
  \tau_{2}(\varphi)
   &\equiv& -\,Tr_{g}R^{N}(\tau(\varphi),d\varphi)d\varphi
  -Tr_{g}\big(\nabla^{\varphi}\nabla^{\varphi}\tau(\varphi)-\nabla^{\varphi}_{\nabla^{M}}\tau(\varphi) \big) = J_{\varphi}(\tau(\varphi)) = 0.
\end{align}

The operator $\tau_{2}(\varphi)$ is called the bitension field of $\varphi$ (see~\cite{C.M.O},~\cite{D.L},~\cite{O.N.D}). It is obvious to see that any harmonic map is biharmonic, therefore it is interesting to construct proper biharmonic maps (non-harmonic biharmonic maps).

 Using the conformal transformation in~\cite{BFS},~\cite{BK} and~\cite{AB}, the authors gives some examples of proper biharmonic maps. In~\cite{BMO} and~\cite{Lu} the authors studied biharmonic maps between warped products where they gave the condition for the
biharmonicity of the inclusion of a Riemannian manifold N into the warped product, also they gave some characterizations of non-harmonic biharmonic maps using the product of harmonic maps and warping metric. The main motivation of this work is to give other methods for the construction of new examples of proper biharmonic maps.

Thus, in this paper, we define the Mus-Gradient metric on tangent bundle $TM$ by a deformation non-conform of Sasaki metric over an n-dimensional Riemannian manifold $(M, g)$ (Definition~\ref{Def_1}). First we investigate the geometry of the Mus-Gradient metric (Theorem~\ref{Th_1} and Theorem~\ref{Th_2}), and we characterize a class of biharmonic maps (Theorem~\ref{Th4.2}, Theorem~\ref{Th4.4}, Theorem~\ref{Th5.2} and Theorem~\ref{Th5.4}). Examples of proper biharmonic maps are constructed
when all of the factors are Euclidean spaces (Example~\ref{Ex4.1}, Example~\ref{Ex4.2}, Example~\ref{Ex5.1} and Example~\ref{Ex5.2}).

\section{Preliminaries}

Let $(M,g)$ be an m-dimensional Riemannian manifold and $(TM,\pi,M)$ be its tangent bundle. Any local chart $(U,x^{i})_{i = \overline{1,m}}$ on $M$ induces a local chart $(\pi^{-1}(U),x^{i},y^{i})_{i = \overline{1,m}}$ on $TM$. Denote by $\Gamma_{ij}^{k}$ the Christoffel symbols of $g$ and by $\nabla$ the Levi-Civita connection of $g$. We have two complementary distributions on $TM$, the vertical distribution $\mathcal{V}$ and the horizontal distribution $\mathcal{H}$, defined by

\[
\mathcal{V}_{(x,u)} = \ker (d\pi_{(x,u)}) = \{a^{i}\frac{\partial}{\partial y^{i}}|_{(x,u)};\; a^{i}\in \mathbb{R}\},
\]
\[
\mathcal{H}_{(x,u)} = \{a^{i}\frac{\partial}{\partial x^{i}}|_{(x,u)}-a^{i}u^{j}\Gamma_{ij}^{k}\frac{\partial}{\partial y^{k}}|_{(x,u)}; a^{i}\in \mathbb{R}\},
\]
where $(x,u)\in TM$,such that $T_{(x,u)}TM = \mathcal{H}_{(x,u)}\oplus \mathcal{V}_{(x,u)}$.
Let $X = X^{i}\frac{\partial}{\partial x^{i}}$ be a local vector field on $M$. The vertical and the horizontal lifts of $X$ are defined by:
\begin{gather}
X^{V} = X^{i}\frac{\partial}{\partial y^{i}},\\
X^{H} = X^{i}\frac{\delta}{\delta x^{i}} = X^{i}\{\frac{\partial}{\partial x^{i}}-y^{j}\Gamma_{ij}^{k}\frac{\partial}{\partial y^{k}}\}. \nonumber
\end{gather}

For consequences, we have ${(\frac{\partial}{\partial x^{i}})}^{H} = \frac{\delta}{\delta x^{i}}$ and ${(\frac{\partial}{\partial x^{i}})}^{V} = \frac{\partial}{\partial y^{i}}$, then $(\frac{\delta}{\delta x^{i}},\frac{\partial}{\partial y^{i}})_{i = \overline{1,m}}$ is a local adapted frame on $TTM$.

\begin{definition}\label{Def_0}
Let $(M,g)$ be a Riemannian manifold. The Sasaki metric $\widehat{g}$ is defined on the tangent bundle $TM$ by
\[
\widehat{g}(X^{H},Y^{H})_{(x,u)} = g_{x}(X,Y), \qquad
\widehat{g}(X^{H},Y^{V})_{(x,u)} = 0, \qquad
\widehat{g}(X^{V},Y^{V})_{(x,u)} = g_{x}(X,Y),
\]
where $X, Y\in\Gamma(TM)$ and $(x,u)\in TM$.
\end{definition}

\begin{theorem}[\cite{Key2}]\label{Theo1.1}
If $\nabla$ (resp. $\widehat{\nabla}$) denote the Levi-Civita connection of $(M,g)$ (resp. $(TM,\widehat{g})$), then for all $X, Y$ and $Z\in\Gamma(TM)$ we have:
\begin{enumerate}[1)]
    \item $(\widehat{\nabla} _{X^{H}}Y^{H})_{(x,u)} = (\nabla _{X}Y)_{(x,u)}^{H}-\frac{1}{2} {(R_{x}(X,Y)u)}^{V}$,

    \item $(\widehat{\nabla} _{X^{H}}Y^{V})_{(x,u)} = (\nabla _{X}Y)_{(x,u)}^{V}+\frac{1}{2}{(R_{x}(u,Y)X)}^{H}$,
    
    \item $(\widehat{\nabla} _{X^{V}}Y^{H})_{(x,u)} = \frac{1}{2}{(R_{x}(u,X)Y)}^{H}$,
    
    \item $\widehat{\nabla} _{X^{V}}Y^{V})_{(x,u)} = 0$,
    
    \item $\widehat{R}_{(x,u)}(X^{V},Y^{V})Z^{V} = 0$,
    
    \item $\widehat{R}_{(x,u)}(X^{V},Y^{V})Z^{H} = [R(X,Y)Z+\frac{1}{4}R(u,X)(R(u,Y)Z)-\frac{1}{4}R(u,Y)(R(u,X)Z)]_{x}^{H}$,
    
    \item $\widehat{R}_{(x,u)}(X^{H},Y^{V})Z^{V} = -[\frac{1}{2}R(Y,Z)X+\frac{1}{4}R(u,Y)(R(u,Z)X)]_{x}^{H}$,
    
    \item $\widehat{R}_{(x,u)}(X^{H},Y^{V})Z^{H} = [\frac{1}{4}R(R(u,Y)Z,X)u+\frac{1}{2}R(X,Z)Y]_{x}^{V}+\frac{1}{2}[(\nabla _{X}R)(u,Y)Z]_{x}^{H}$,
    
    \item 
    $\begin{aligned}[t]
    \widehat{R}_{(x,u)}(X^{H},Y^{H})Z^{V} 
    ={ }&{ }[R(X,Y)Z+\frac{1}{4}R(R(u,Z)Y,X)u-\frac{1}{4}R(R(u,Z)X,Y)u]_{x}^{V}\\
    &{ }+\frac{1}{2}[(\nabla _{X}R)(u,Z)Y-(\nabla _{Y}R)(u,Z)X]_{x}^{H},
    \end{aligned}$

    \item 
    $\begin{aligned}[t]
    \widehat{R}_{(x,u)}(X^{H},Y^{H})Z^{H}
    ={ }&{ } \frac{1}{2}[(\nabla _{Z}R)(X,Y)u]^{V}_{x}+[R(X,Y)Z+\frac{1}{4}R(u,R(Z,Y)u)X\\
    &{ }+\frac{1}{4}R(u,R(X,Z)u)Y+\frac{1}{2}R(u,R(X,Y)u)Z]^{H}_{x},
    \end{aligned}$
\end{enumerate}
where $(x,u)\in TM$ such that $\pi(u) = x$ and $R$ (resp $\widehat{R}$) denote the curvature tensor of $(M,g)$ (resp $(M,\widehat{g})$).
\end{theorem}

\begin{remark}
From Theorem~\ref{Theo1.1} we conclude that $TM$ is a rigid fiber bundle, i.e $(TM,\widehat{g})$ is flat if and only if $(M,g)$ is flat.
\end{remark}
\begin {corollary}\label{Cor.0.1} 
If $(M,g)$ is flat ($R = 0$), then we obtain
\begin{enumerate}[1)]
    \item $(\widehat{\nabla} _{X^{H}}Y^{H})_{(x,u)} = (\nabla _{X}Y)_{(x,u)}^{H}$,
    
    \item $(\widehat{\nabla} _{X^{H}}Y^{V})_{(x,u)} = (\nabla _{X}Y)_{(x,u)}^{V}$,
    
    \item $(\widehat{\nabla} _{X^{V}}Y^{H})_{(x,u)} = 0$,
    
    \item $(\widehat{\nabla} _{X^{V}}Y^{V})_{(x,u)} = 0$,
    
    \item $\widehat{R} = 0$.
\end{enumerate}
\end{corollary}

The geometry of the tangent bundle $TM$ equipped withe Sasaki metric has been studied by many authors Sasaki~\cite{Key2}, K.Yano and S. Ishihara~\cite{Key7}, P.Dombrowski~\cite{Key3}, A. A. Salimov, and Gezer~\cite{Key19},~\cite{Key16},~\cite{Key15},~\cite{Key17} and others
(see~\cite{Aba},~\cite{Key5},~\cite{Key9},~\cite{Key8},~\cite{Key4},~\cite{Z.D}).

\begin{definition}[\cite{Z.D}]\label{Def.0.2}
Let $(M,g)$ be a Riemannian manifold. The Mus-Sasaki metric $\widetilde{g}$ is defined on the tangent bundle $TM$ by
\begin{enumerate}
\item $\widetilde{g}(X^{H},Y^{H})_{(x,u)} = g_{x}(X,Y)$,
\item $\widetilde{g}(X^{H},Y^{V})_{(x,u)} = 0$,
\item $\widetilde{g}(X^{V},Y^{V})_{(x,u)} = f(x)g_{x}(X,Y)$,
\end{enumerate}
where $X, Y\in\Gamma(TM)$ and $(x,u)\in TM$.
\end{definition}

\begin {lemma}\label{Lem.0.2}
Let $(M,g)$ be a flat Riemannian manifold. Then we have
\begin{enumerate}[1)]
    \item $(\widetilde{\nabla} _{X^{H}}Y^{H})_{(x,u)} = (\nabla _{X}Y)_{(x,u)}^{H}$,

    \item $(\widetilde{\nabla} _{X^{H}}Y^{V})_{(x,u)} = (\nabla _{X}Y)_{(x,u)}^{V}+\frac{1}{2}X(f)Y^{V}$,

    \item $\widetilde{\nabla} _{X^{V}}Y^{H})_{(x,u)} = \frac{1}{2}Y(f)X^{V}$,

    \item $(\widetilde{\nabla} _{X^{V}}Y^{V})_{(x,u)} = -\frac{1}{2}g(X,Y)\big{(grad(f)\big)}^{H}$,
    
    \item $\widetilde{R}(X^{H},Y^{H})Y^{H} = 0$,
    
    \item $\widetilde{R}(X^{H},Y^{V})Y^{V} = -\frac{1}{2}\|Y\|^{2}\big{(\nabla_{X}grad(f)\big)}^{H}+ \frac{1}{4f}\|Y\|X(f)\big{(grad(f)\big)}^{H}$,
\end{enumerate}
where $\widetilde{\nabla}$ (resp $\widetilde{R}$) denote the Levi-Civita connection (resp curvature tensor) of $\widetilde{g}$.
\end{lemma}

For more detail on geometry of Mus-Sasaki metric see~\cite{L.D.Z},~\cite{Z.D}.

\section{Mus-Gradient metric}

\begin{definition}\label{Def_1}
Let $(M,g)$ be a Riemannian manifold and $f:M \rightarrow ]0,+\infty[$ be a function such that $f\in C^{\infty}(M)$. On the tangent bundle $TM$, we define a Mus-Gradient metric noted $g^{f}$ by:
\begin{enumerate}
\item $g^{f}(X^{H},Y^{H})_{(x,u)} = g_{x}(X,Y)$,
\item $g^{f}(X^{H},Y^{V})_{(x,u)} = 0$,
\item $g^{f}(X^{V},Y^{V})_{(x,u)} = g_{x}(X,Y)+X_{x}(f)Y_{x}(f)$,
\end{enumerate}
where $X, Y\in\Gamma(TM)$, $(x,u)\in TM$.
\end{definition}

\begin{remark}\label{re_1}
\begin{enumerate}[1)]
\item If $f$ is constant then $g^{f} = g^{S}$ is the Sasaki metric~\cite{Key4},
\item $g^{f}(X^{H},(grad{(f)}^{H}) = g(X,grad(f) = X(f)$,
\item $g^{f}(X^{V},(grad{(f)}^{V}) = (1+\|\operatorname{grad}f\|^{2})X(f)$,
\item $g^{f}(X^{V},Y^{V})-g^{f}(X^{H},Y^{H}) = X(f)Y(f)$,
\end{enumerate}
where $X, Y\in\Gamma(TM)$.
\end{remark}

In the following, we consider $\alpha = 1+\|\operatorname{grad}f\|^{2} = 1+g(grad\;f,grad\;f)$.

\begin{lemma}\label{Lem_1}
Let $(M,g)$ be a Riemannian manifold. for all $X, Y$ and $Z\in\Gamma(TM)$, we have
\begin{enumerate}[1)]
    \item $X^{V}\big(g^{f}(Y^{V},Z^{V})\big) = 0$,

    \item 
    $\begin{aligned}[t]
    X^{H}\big(g^{f}(Y^{V},Z^{V})\big) 
    ={ }&{ }g^{f}({(\nabla_{X}Y)}^{V},Z^{V}) + g^{f}({(Y^{V},\nabla_{X}Z)}^{V})\\
    &+ Y(f)g^{f}((\nabla_{X}grad{(f)}^{H},Z^{H})\\
    &+ g(Y,\nabla_{X}grad(f)g^{f}((grad{(f)}^{H},Z^{H}),
    \end{aligned}$

    \item
    $\begin{aligned}[t]
    X^{H}\big(g^{f}(Y^{V},Z^{V})\big)
    ={ }&{ }g^{f}({(\nabla_{X}Y)}^{V},Z^{V}) + g^{f}({(Y^{V},\nabla_{X}Z)}^{V})\\
    &+ Y(f)g^{f}((\nabla_{X} grad{(f)}^{V},Z^{V}) + \frac{1}{\alpha}\big[g(Y,\nabla_{X}grad(f)\\
    &-\frac{1}{2}X(\alpha) Y(f)\big]g^{f}((grad{(f)}^{V},Z^{V}).
    \end{aligned}$
\end{enumerate}
\end{lemma}

\begin{proof}
\quad The statement is a direct consequence of Definition~\ref{Def_1} and Remark~\ref{re_1}.
\end{proof}

\begin{lemma}\label{Lem_1.1}
Let $(M,g)$ be a Riemannian manifold. If $\nabla$ (resp. $\nabla^{f}$) denote the Levi-Civita connection of $(M,g)$ (resp. $(TM,g^{f})$), then for all $X, Y$ and $Z\in\Gamma(TM)$ we have:
\begin{enumerate}[1)]
    \item $g^{f}(\nabla^{f}_{X^{H}}Y^{H},Z^{H}) = g^{f}({(\nabla_{X}Y)}^{H},Z^{H})$,

    \item $g^{f}(\nabla^{f}_{X^{H}}Y^{H},Z^{V}) = -\frac{1}{2}g^{f}({(R(X,Y)u)}^{V},Z^{V})$,

    \item $g^{f}(\nabla^{f}_{X^{H}}Y^{V},Z^{H}) = \frac{1}{2}g^{f}({(R(u,Y)X)}^{H}+Y(f)(R{(u,grad(f)X)}^{H},Z^{H})$,

    \item 
    $\begin{aligned}[t]
    g^{f}(\nabla^{f}_{X^{H}}Y^{V},Z^{V}) 
    ={ }&{ }g^{f}({(\nabla_{X}Y)}^{V},Z^{V})+\frac{1}{2}Y(f)g^{f}((\nabla_{X} grad{(f)}^{V},Z^{V})\\
    &+\frac{1}{2\alpha}\,\big[g(Y,\nabla_{X}grad(f)-\frac{1}{2}X(\alpha)Y(f)\big]g^{f}((grad{(f)}^{V},Z^{V}),
    \end{aligned}$
    
    \item $g^{f}(\nabla^{f}_{X^{V}}Y^{H},Z^{H}) = \frac{1}{2}g^{f}({(R(u,X)Y)}^{H}+X(f)(R{(u,grad(f)Y)}^{H},Z^{H})$,

    \item
    $\begin{aligned}[t]
    g^{f}(\nabla^{f}_{X^{V}}Y^{H},Z^{V})
    ={ }&{ }\frac{1}{2}X(f)g^{f}((\nabla_{Y} grad{(f)}^{V},Z^{V}) \\
    &+\frac{1}{2\alpha}\big[g(X,\nabla_{Y}grad(f) - \frac{1}{2}Y(\alpha)X(f)\big]g^{f}((grad{(f)}^{V},Z^{V}),
    \end{aligned}$

    \item $g^{f}(\nabla^{f}_{X^{V}}Y^{V},Z^{H}) = -\frac{1}{2}g^{f}\big(X(f)(\nabla_{Y} grad{(f)}^{H}+Y(f)(\nabla_{X} grad{(f)}^{H},Z^{H}\big)$,

    \item $g^{f}(\nabla^{f}_{X^{V}}Y^{V},Z^{V}) = 0$.
\end{enumerate}
\end{lemma}

\begin{proof}
The proof of Lemma~\ref{Lem_1.1} follows directly from Koszul formula, Definition~\ref{Def_1} and Lemma~\ref{Lem_1}.
\end{proof}

As a direct consequence of Lemma~\ref{Lem_1.1}, we get the following theorem.
\begin{theorem}\label{Th_1}
Let $(M,g)$ be a Riemannian manifold. If $\nabla$ (resp. $\nabla^{f}$) denote the Levi-Civita connection of $(M,g)$ (resp. $(TM,g^{f})$), then we have
\begin{enumerate}[1)]
    \item $(\nabla^{f}_{X^{H}}Y^{H})_{p} = {(\nabla_{X}Y)}^{H}_{p}-\frac{1}{2}{(R_{x}(X,Y)u)}^{V}$,

    \item 
    $\begin{aligned}[t]
    (\nabla^{f}_{X^{H}}Y^{V})_{p}
    =& \frac{1}{2}{(R_{x}(u,Y)X)}^{H}+\frac{1}{2}Y_{x}(f)(R_{x}{(u,grad(f)X)}^{H}+\frac{1}{2}Y_{x}(f)(\nabla_{X} grad{(f)}^{V}_{p}\\
    &+{(\nabla_{X}Y)}^{V}_{p}+\frac{1}{2\alpha}\big[g_{x}(Y,\nabla_{X}grad(f)-\frac{1}{2}X_{x}(\alpha)Y_{x}(f)\big](grad{(f)}^{V}_{p},
    \end{aligned}$

    \item 
    $\begin{aligned}[t]
    (\nabla^{f}_{X^{V}}Y^{H})_{p}
    =& \frac{1}{2}{(R_{x}(u,X)Y)}^{H}+\frac{1}{2}X_{x}(f)(R_{x}{(u,grad(f)Y)}^{H}+\frac{1}{2}X_{x}(f)(\nabla_{Y} grad{(f)}^{V}_{p}\\
    &+\frac{1}{2\alpha}\big[g_{x}(X,\nabla_{Y}grad(f)-\frac{1}{2}Y_{x}(\alpha)X_{x}(f)\big](grad{(f)}^{V}_{p},
    \end{aligned}$

    \item $(\nabla^{f}_{X^{V}}Y^{V})_{p} = -\frac{1}{2}X_{x}(f)(\nabla_{Y} grad{(f)}^{H}_{p}-\frac{1}{2}Y_{x}(f)(\nabla_{X} grad{(f)}^{H}_{p}$,
\end{enumerate}
for all vector fields $X,Y\in \Gamma(TM)$ and $p = (x,u)\in TM$, where $R$ denote the curvature tensor of $(M,g)$.
\end{theorem}

\begin{corollary}\label{Cor.1-2}
Let $(M,g)$ be a flat manifold. If $\nabla$ (resp. $\nabla^{f}$) denote the Levi-Civita connection of $(M,g)$ (resp. $(TM,g^{f})$), then we have:
\begin{enumerate}[1)]
    \item $(\nabla^{f}_{X^{H}}Y^{H})_{p} = {(\nabla_{X}Y)}^{H}_{p}$,

    \item 
    $\begin{aligned}[t]
    (\nabla^{f}_{X^{H}}Y^{V})_{p} 
    =& \frac{1}{2}Y_{x}(f)(\nabla_{X} grad{(f)}^{V}_{p}+{(\nabla_{X}Y)}^{V}_{p}+\frac{1}{2\alpha}\big[g_{x}(Y,\nabla_{X}grad(f)\\
    &-\frac{1}{2}X_{x}(\alpha)Y_{x}(f)\big](grad{(f)}^{V}_{p},
    \end{aligned}$

    \item 
    $\begin{aligned}[t]
    (\nabla^{f}_{X^{V}}Y^{H})_{p}
    =& \frac{1}{2}X_{x}(f)(\nabla_{Y} grad{(f)}^{V}_{p} \\
    &+\frac{1}{2\alpha}\Big[g_{x}\big(X,\nabla_{Y}grad(f)\big)-\frac{1}{2}Y_{x}(\alpha)X_{x}(f)\Big](grad{(f)}^{V}_{p},
    \end{aligned}$

    \item $(\nabla^{f}_{X^{V}}Y^{V})_{p} = -\frac{1}{2}X_{x}(f)(\nabla_{Y} grad{(f)}^{H}_{p}-\frac{1}{2}Y_{x}(f)(\nabla_{X} grad{(f)}^{H}_{p}$,
 \end{enumerate}
for all vector fields $X,Y\in \Gamma(TM)$ and $p = (x,u)\in TM$, where $R$ denote the curvature tensor of $(M,g)$.
\end{corollary}

From Theorem~\ref{Th_1} we obtain the following theorem

\begin{theorem}\label{Th_2}
Let $(M,g)$ be a Riemannian manifold and $(TM,g^{f})$) its tangent bundle equipped with the Mus-Gradient metric. If $R$ (resp $R^{f}$ denote the Riemannian curvature tensor of $(M,g)$ (resp $(TM,g^{f})$), then we have the following formulas:
\begin{align*}
R^{f}_{p}&(X^{H},Y^{H})Y^{H}\\
=& {(R_{x}(X,Y)Y)}^{H}+\frac{3}{4}{(R_{x}(u,R(X,Y)u)Y)}^{H}\\
&+\frac{3}{4}g_{x}(R(X,Y)u,grad(f)(R_{x}{(u,grad(f)Y)}^{H}\\
&+\frac{1}{2}{((\nabla _{Y}R)(X,Y)u)}^{V}_{p}+\frac{3}{4}g_{x}(R(X,Y)u,grad(f) (\nabla _{Y}grad{(f)}^{V}_{p}\\
&+\frac{3}{8\alpha}\big[2g_{x}(R(X,Y)u,\nabla _{Y}grad(f) -Y_{x}(\alpha)g_{x}(R(X,Y)u,grad(f))\big](grad{(f)}^{V}_{p},
\\
R^{f}_{p}&(X^{H},Y^{V})Y^{V}\\
=& -Y_{x}(f)(\nabla_{X}\nabla_{Y}grad{(f)}^{H}_{p}-\frac{1}{4}Y_{x}(f) \big(R_{x}{(u,Y+grad(f)R(u,Y)X\big)}^{H}\\
&-\frac{1}{2} Y_{x}(f)(R_{x}{(Y,grad(f)X)}^{H}-\frac{1}{4}Y_{x}(f)(R_{x}(u,Y+Y(f)grad(f)R{(u,grad(f)X)}^{H}\\
&+\frac{1}{4}Y^{2}_{x}(f)(\nabla_{(\nabla_{X}grad(f)}grad{(f)}^{H}_{p}+Y_{x}(f) (\nabla_{(\nabla_{X}Y)}grad{(f)}^{H}_{p}\\
&+\frac{Y_{x}(f)}{8\alpha}g_{x}\big(Y,\nabla_{X}grad(f)-\frac{1}{2}X(\alpha)grad(f)\big)\big{(grad(\alpha)\big)}^{H}_{p} \\
&+\big[\frac{1}{8\alpha}X_{x}(\alpha)Y_{x}(f)-\frac{1+3\alpha}{4\alpha}g_{x}(Y,\nabla_{X}grad(f)\big] (\nabla _{Y}grad{(f)}^{H}_{p}\\
&+\frac{1}{2}Y_{x}(f)(R_{x}{(X,\nabla_{Y}grad(f)u)}^{V}
-\frac{1}{4}Y_{x}(f)(\nabla_{(R(u,Y+Y(f)grad(f)X)}grad{(f)}^{V}_{p} \\
&+\frac{1}{8\alpha}g_{x}\big(R(u,Y+Y(f)grad(f)X,Y(f)grad\alpha-\nabla_{Y}grad(f)\big)(grad{(f)}^{V}_{p},
\\
-R^{f}_{p}&(X^{V},Y^{H})Y^{H}\\
=& \frac{1}{2}{((\nabla_{Y}R)_{x}(u,X)Y)}^{H}+\frac{1}{2}X_{x}(f)((\nabla_{Y}R)_{x}{(u,grad(f)Y)}^{H} \\
&+\frac{3}{4}X_{x}(f)(R_{x}{(u,\nabla_{Y}grad(f)Y)}^{H}+\frac{3}{4}g_{x}(X,\nabla_{Y}grad(f)(R_{x}{(u,grad(f)Y)}^{H} \\
&+\frac{1}{2}X_{x}(f)\big[\nabla_{Y}\nabla_{Y}grad(f)-\nabla_{(\nabla_{Y}Y)}grad(f)\big]^{V}_{p}\\
&-\frac{1}{4}\big(R_{x}{(Y,R(u,(X+X(f)grad(f))Y)u\big)}^{V} \\
&+\big[\frac{3\alpha+1}{4\alpha}g_{x}(X,\nabla_{Y}grad(f)-\frac{1}{8\alpha}X_{x}(f)Y_{x}(\alpha)\big] (\nabla _{Y}grad{(f)}^{V}_{p}\\
&+\big[\frac{1}{4\alpha}X_{x}(f)\|\nabla_{Y}grad(f)\|^{2}+\frac{1}{2\alpha}g_{x}\big(\nabla_{Y}\nabla_{Y}grad(f)-\nabla_{(\nabla_{Y}Y)}\,grad(f),X\big) \\
&-\frac{1}{4\alpha}X_{x}(f)g_{x}(Y,\nabla_{Y}grad(\alpha)) -\frac{1}{8\alpha^{2}}X_{x}(f){(Y_{x}(\alpha))}^{2} \\
&-\frac{3\alpha+2}{8\alpha^{2}}Y_{x}(\alpha)g_{x}(X,\nabla_{Y}grad(f)\big](grad{(f)}^{V}_{p},
\\
R^{f}_{p}&(X^{V},Y^{V})Y^{V}\\
=& \frac{1}{4}(R_{x}(u,Y){(\nabla _{\big(X(f)Y+Y(f)X\big)}grad(f))}^{H}-\frac{1}{2}Y_{x}(f)(R_{x}(u,X){(\nabla_{Y}grad(f))}^{H} \\
&+\frac{1}{4}Y_{x}(f)\big(R_{x}(u,grad(f){(\nabla_{\big(Y(f)X-X(f)Y\big)}grad(f)\big)}^{H} \\
&+\frac{1}{8\alpha} g_{x} \left( \nabla_{\left( X(f)Y-Y(f)X \right)} grad(f), \left( Y(f) grad(\alpha) + 2\nabla _Y grad(f) \right) grad{(f)}^{V}_{p} \right) \\
&+\frac{1}{4} Y_{x}(f) \left( \nabla_{\left( \nabla _{\left( Y(f)X-X(f)Y \right)} grad(f) \right)} grad(f) \right)^{V}_{p},
\end{align*}
for all $p = (x,u)\in TM$ and $,X,Y,Z\in \Gamma(TM)$.
\end{theorem}

Using Theorem~\ref{Th_2}, we obtain the following corollaries:

\begin{corollary}\label{Cor_2.1}
If $(M,g)$ is a flat Riemannian manifold, then we have:
\begin{align*}
R^{f}_{p}(X^{H},Y^{H})Y^{H} 
=& 0,
\\
R^{f}_{p}(X^{H},Y^{V})Y^{V}
=& -Y_{x}(f)(\nabla_{X}\nabla_{Y}grad{(f)}^{H}_{p}\\
&+\frac{1}{4}Y^{2}_{x}(f)\big{(\nabla_{\nabla_{X}grad(f)}grad(f)\big)}^{H}_{p}+Y_{x}(f) (\nabla_{(\nabla_{X}Y)}grad{(f)}^{H}_{p}\\
&+\frac{Y_{x}(f)}{8\alpha}g_{x}\big(Y,\nabla_{X}grad(f)-\frac{1}{2}X(\alpha)grad(f)\big){(grad(\alpha))}^{H}_{p} \\
&+\big[\frac{1}{8\alpha}X_{x}(\alpha)Y_{x}(f)-\frac{1+3\alpha}{4\alpha}g_{x}(Y,\nabla_{X}grad(f)\big] (\nabla _{Y}grad{(f)}^{H}_{p},
\\
-R^{f}_{p}(X^{V},Y^{H})Y^{H}
=& \frac{1}{2}X_{x}(f)\big[\nabla_{Y}\nabla_{Y}grad(f)-\nabla_{(\nabla_{Y}Y)}grad(f)\big]^{V}_{p} \\
&+\big[\frac{3\alpha+1}{4\alpha}g_{x}(X,\nabla_{Y}grad(f)-\frac{1}{8\alpha}X_{x}(f)Y_{x}(\alpha)\big] (\nabla _{Y}grad{(f)}^{V}_{p}\\
&+\big[\frac{1}{4\alpha}X_{x}(f)\|\nabla_{Y}grad(f)\|^{2}\\
&+\frac{1}{2\alpha}g_{x}\big(\nabla_{Y}\nabla_{Y}grad(f)-\nabla_{(\nabla_{Y}Y)}\,grad(f),X\big) \\
&-\frac{1}{4\alpha}X_{x}(f)g_{x}(Y,\nabla_{Y}grad(\alpha)) -\frac{1}{8\alpha^{2}}X_{x}(f){(Y_{x}(\alpha))}^{2} \\
&-\frac{3\alpha+2}{8\alpha^{2}}Y_{x}(\alpha)g_{x}(X,\nabla_{Y}grad(f)\big](grad{(f)}^{V}_{p},
\\
R^{f}_{p}(X^{V},Y^{V})Y^{V}\\
={\ \ \ }&{\!\!\!\!\!} \frac{1}{8\alpha}g_{x}\big(\nabla_{\big(X(f)Y-Y(f)X\big)}grad(f),(Y(f)grad(\alpha)+2\nabla _{Y}grad(f)\big)(grad{(f)}^{V}_{p}\\
&{\!\!\!\!\!}+\frac{1}{4}Y_{x}(f) \left( \nabla_{\left( \nabla _{\left( Y(f)X-X(f)Y \right)} grad(f) \right)} gradf \right)^{V}_{p},
\end{align*}
for all $p = (x,u)\in TM$ and $,X,Y,Z\in \Gamma(TM)$.
\end{corollary}

\begin{corollary}\label{Cor_2.2} 
If $(M,g)$ is a flat Riemannian manifold then $(TM,g^{f})$ is flat if and only if $f$ is constant.
\end{corollary}

\section{Biharmonic maps of Mus-Gradient metric}

\begin{lemma}\label{Le4.1}
Let $(M,g)$ be a Riemannian manifold and $(TM,g^{S})$ be its tangent bundle equipped with the Sasaki metric then the canonical projection $\pi:(TM,g^{S})\mapsto (M,g)$ is harmonic.
\end{lemma}

\begin{theorem}\label{Th4.1}
Let $(M,g)$ be a Riemannian manifold, $f:M \rightarrow ]0,+\infty[$ such $f\in C^{\infty}(M)$ and $(TM,g^{f})$ be its tangent bundle equipped with the Mus-Gradient metric. The tension field of the canonical projection $\pi_{\alpha}:(TM,g^{f})\mapsto (M,g)$ is given by
\begin{equation}\label{Eq4.3}
    \tau(\pi_{\alpha}) = \frac{1}{\alpha}\nabla^{}_{gradf} gradf = \frac{grad\alpha}{2\alpha},
\end{equation}
and $\pi_{\alpha}$ is harmonic if and only if $\|grad f\| = Const$.
\end{theorem}

\begin{proof}
If $(E_{1},\dotsc,E_{m})$ is orthonormal frame on $(M,g)$ such that $E_{1} = \frac{grad f}{\|grad f\|}$, then
\[
\left( E_{1}^{H},\dotsc,E^{H}_{m},\frac{1}{\sqrt{\alpha}}E_{1}^{V},\dotsc,E_{m}^{V}\right)
\]
is an orthonormal frame on $(TM,g^{f})$.  From formula (\ref{Eq0.2}) and Theorem~\ref{Th_1} we obtain:
\begin{align*}
  \tau(\pi_{\alpha})
  &= \sum_{i = 1}^{m}\nabla_{E_{i}}E_{i} - \sum_{i = 1}^{m}d\pi_{\alpha}(\nabla^{f}_{E^{H}_{i}}E^{H}_{i})-\sum_{i = 2}^{m}d\pi_{\alpha}(\nabla^{f}_{E^{V}_{i}}E^{V}_{i})-\frac{1}{\alpha}d\pi_{\alpha}(\nabla^{f}_{E^{V}_{1}}E^{V}_{1})\\
   &= \sum_{i = 2}^{m}E_{i}(f)\nabla_{E_{i}}gradf+\frac{1}{\alpha}E_{1}(f)\nabla_{E_{1}}gradf\\
   &= \nabla_{\sum_{i = 1}^{m}E_{i}(f)E_{i}}gradf+\frac{1-\alpha}{\alpha}E_{1}(f)\nabla_{E_{1}}gradf\\
   &= \nabla_{gradf}gradf+\frac{1-\alpha}{\alpha}E_{1}(f)\nabla_{E_{1}}gradf\\
   &= \nabla_{gradf}gradf+\frac{1-\alpha}{\alpha}\|gradf\|\nabla_{\frac{gradf}{\|gradf\|}}gradf\\
   &= \nabla_{gradf}gradf+\frac{1-\alpha}{\alpha}\nabla_{gradf}gradf\\
    &= \frac{1}{\alpha}\nabla_{gradf}gradf\\
    &= \frac{grad\alpha}{2\alpha}.
\qedhere
\end{align*}
\end{proof}

\begin{theorem}\label{Th4.2}
Let $(M,g)$ be a Riemannian manifold, $f$ a smooth positive function on $M$ such that $gradf\neq0$ at any point on $M$, and $(TM,g^{f})$ its tangent bundle equipped with the Mus-Gradient metric. Then the canonical projection $\pi_{\alpha}:(TM,g^{f})\longrightarrow (M,g)$ is biharmonic if and only if
\begin{align}\label{Eq4.4}
   \operatorname{Ricci}(grad\ln\alpha) +\frac{1}{2}grad(\Delta \ln\alpha)
       +\frac{1}{8}grad(\|grad\ln\alpha\|^2) = 0.
\end{align}
\end{theorem}

\begin{proof}
Let $(E_{1},\dots,E_{m})$ be an orthonormal frame on $(M,g)$ such that $E_{1} = \frac{grad f}{\|grad f\|}$, so that
 $(E_{1}^{H},\dots,E^{H}_{m},\frac{1}{\sqrt{\alpha}}E_{1}^{V},\dots,E_{m}^{V})$
is an orthonormal frame on $(TM,g^{f})$. Since $d\pi_{\alpha}(X^H) = X$ and $d\pi_{\alpha}(X^V) = 0$ for all $X\in\Gamma(TM)$, we obtain
\begin{align*}
-\operatorname{Tr}_{g^f}R^{}(\tau(\pi_{\alpha}),d\pi_{\alpha})d\pi_{\alpha} = -\frac{1}{2\alpha}\sum_{i = 1}^m R(grad\alpha,E_i)E_i,
\end{align*}
be the definition of Ricci tensor of $(M,g)$ with $\frac{1}{\alpha}grad\alpha = grad\ln\alpha$, we get
\begin{align*}
-\operatorname{Tr}_{g^f}R^{}(\tau(\pi_{\alpha}),d\pi_{\alpha})d\pi_{\alpha} = -\frac{1}{2}\operatorname{Ricci}(grad\ln\alpha).
\end{align*}
 The term $-\operatorname{Tr}_{g^f}\nabla^{\pi_{\alpha}}\nabla^{\pi_{\alpha}}\tau(\pi_{\alpha})$ is given by
\[-\operatorname{Tr}_{g^f}\nabla^{\pi_{\alpha}}\nabla^{\pi_{\alpha}}\tau(\pi_{\alpha}) = -\frac{1}{2}\sum_{i = 1}^m \nabla^{\pi_{\alpha}}_{E_i^H}\nabla^{\pi_{\alpha}}_{E_i^H}grad\ln\alpha = -\frac{1}{2}\sum_{i = 1}^m \nabla^{}_{E_i}\nabla^{}_{E_i}grad\ln\alpha.
\]
For the last term $\operatorname{Tr}_{g^f}\nabla^{\pi_{\alpha}}_{\nabla^{f}}\tau(\pi_{\alpha})$ we have
\begin{equation}\label{Eq4.5}
  \operatorname{Tr}_{g^f}\nabla^{\pi_{\alpha}}_{\nabla^{f}}\tau(\pi_{\alpha})
   = \frac{1}{2}\sum_{i = 1}^m \nabla^{\pi_{\alpha}}_{\nabla^{f}_{E_i^H}E_i^H}grad\ln\alpha
       +\frac{1}{2\alpha}\nabla^{\pi_{\alpha}}_{\nabla^{f}_{E_1^V}E_1^V}grad\ln\alpha
 +\frac{1}{2}\sum_{i = 2}^m \nabla^{\pi_{\alpha}}_{\nabla^{f}_{E_i^V}E_i^V}grad\ln\alpha.
\end{equation}
According to Theorem~\ref{Th_1}, we obtain
\begin{align*}
 (\nabla^{f}_{E_i^{H}}E_i^{H})_{(x,u)} &= {(\nabla_{E_i}E_i)}^{H}_{x}-\frac{1}{2}{(R_{x}(E_i,E_i)u)}^{V} = {(\nabla_{E_i}E_i)}^{H}_{x}, \hbox{for all $1\leq i \leq m$}\\
 (\nabla^{f}_{E_1^{V}}E_1^{V})_{(x,u)} &= -E_1(f)_{x}(\nabla_{E_1} grad{(f)}^{H}_{x} = -(\nabla_{gradf} grad{(f)}^{H}_{x},\\
 (\nabla^{f}_{E_i^{V}}E_i^{V})_{(x,u)} &= -E_i(f)_{x}(\nabla_{E_i} grad{(f)}^{H}_{x} = 0, \hbox{for all $2\leq i \leq m$},
\end{align*}
where $(x,u)\in TM$. Substituting the last formulas in (\ref{Eq4.5}) and using the identity 
\[
\nabla_{gradf} gradf = \frac{1}{2}grad(\|gradf\|^2),
\]
we find that
\begin{align*}
  \operatorname{Tr}_{g^f}\nabla^{\pi_{\alpha}}_{\nabla^{f}}\tau(\pi_{\alpha})
   &= \frac{1}{2}\sum_{i = 1}^m \nabla^{}_{\nabla^{}_{E_i}E_i}grad\ln\alpha
       -\frac{1}{4\alpha}\nabla^{}_{grad(\|gradf\|^2)}grad\ln\alpha\\
   &= \frac{1}{2}\sum_{i = 1}^m \nabla^{}_{\nabla^{}_{E_i}E_i}grad\ln\alpha
       -\frac{1}{4}\nabla^{}_{grad\ln\alpha}grad\ln\alpha\\
   &= \frac{1}{2}\sum_{i = 1}^m \nabla^{}_{\nabla^{}_{E_i}E_i}grad\ln\alpha
       -\frac{1}{8}grad(\|grad\ln\alpha\|^2).
\end{align*}
Note that the canonical projection $\pi_{\alpha}$ is biharmonic if and only if
\[
\tau_{2}(\pi_{\alpha})
   = -\operatorname{Tr}_{g^f}R^{}(\tau(\pi_{\alpha}),d\pi_{\alpha})d\pi_{\alpha}
    -\operatorname{Tr}_{g^f}\big(\nabla^{\pi_{\alpha}}\nabla^{\pi_{\alpha}}\tau(\pi_{\alpha})-\nabla^{\pi_{\alpha}}_{\nabla^{}}\tau(\pi_{\alpha}) \big) = 0,
\]
is equivalent to following equation
\begin{align*}
  0 =& -\frac{1}{2}\operatorname{Ricci}(grad\ln\alpha) -\frac{1}{2}\sum_{i = 1}^m \nabla^{}_{E_i}\nabla^{}_{E_i}grad\ln\alpha\\
   & +\frac{1}{2}\sum_{i = 1}^m \nabla^{}_{\nabla^{}_{E_i}E_i}grad\ln\alpha
       -\frac{1}{8}grad(\|grad\ln\alpha\|^2).
\end{align*}
The Theorem~\ref{Th4.2} follows from the last equation
and the following
 \[
\operatorname{Tr}_g\nabla^{2}grad\ln\alpha = \operatorname{Ricci}(grad\ln\alpha)+grad(\Delta \ln\alpha).
\qedhere
\]
\end{proof}

\begin{corollary}\label{co2.}
Let $(M,g)$ be an Einstein manifold (i.e. $\operatorname{Ricci}(X) = \lambda X$ for all $X$ in $\Gamma(TM)$ where $\lambda\in\mathbb{R}$), $f$ a smooth positive function on $M$ such that $gradf\neq0$ at any point on $M$, and $(TM,g^{f})$ its tangent bundle equipped with the Mus-Gradient metric. Then the canonical projection $\pi_{\alpha}:(TM,g^{f})\longrightarrow (M,g)$ is biharmonic if and only if
\begin{align}\label{Eq4.6}
   \lambda\ln\alpha +\frac{1}{2}\Delta \ln\alpha
       +\frac{1}{8}\|grad\ln\alpha\|^2 = Const.
\end{align}
\end{corollary}

According to Corollary~\ref{co2.}, we have the following example.

\begin{example}\label{Ex4.1}
Let $M = (1,\infty)\times\mathbb{R}^{m-1}$ equipped with the Riemannian metric given by $g = dt^2 +\sum_{i = 1}^{m-1}dx_i^2$. Note that, $(M,g)$ is an Einstein manifold with $\lambda = 0$.
We set
\[
f(t,x) = \int \sqrt{t^4 -1}dt, \quad \hbox{for all $(t,x)$ in $M$}.
\]
We have, $grad (f) = \sqrt{t^4 -1}\,\partial_t$, $\|gradf\| = \sqrt{t^4 -1}$, $\alpha(t,x) = t^4$, for all $(t,x)$ in $M$, so that
\[
\Delta \ln\alpha = -\frac{4}{t^2} \quad \hbox{and} \quad \|grad\ln\alpha\|^2 = \frac{16}{t^2}.
\]
Using the Corollary~\ref{co2.}, the canonical projection $\pi_{\alpha}:(TM,g^{f})\longrightarrow (M,g)$ is proper biharmonic. Here, the tension field of $\pi_{\alpha}$ and the Mus-Gradient metric $g^f$ are given by
\[
\tau(\pi_{\alpha}) = \frac{2}{t}\,\partial_t \quad\hbox{and}\quad
g^f = dt^2 +\sum_{i = 1}^{m-1}dx_i^2 +t^4 dy_1^2 +\sum_{i = 1}^{m-1}dy_i^{2}.
\]
\end{example}

\begin{theorem}\label{Th4.3}
Let $(M,g)$ be a flat Riemannian manifold, $f:M \rightarrow ]0,+\infty[$ be a smooth function on $M$ and $(TM,g^{f})$) its tangent bundle equipped with the Mus-Gradient metric. Then the tension field of the identity $Id_{\alpha}:(TM,g^{f})\mapsto (TM,\widehat{g})$ is given by
\begin{equation}\label{Eq.4.1.1}
    \tau(Id_{\alpha}) = \frac{\|grad(f)\|}{1+\|grad(f)\|^{2}}\big{(grad(\|grad(f)\|)\big)}^{H},
\end{equation}
and $Id_{\alpha}$ is harmonic if and only if $\|grad (f)\| = Const$.
\end{theorem}

\begin{proof}
If $(E_{1},\dotsc,E_{m})$ is orthonormal frame on $(M,g)$ such that $E_{1} = \frac{grad f}{\|grad f\|}$, then
\[
\left( E_{1}^{H},\dotsc ,E^{H}_{m},\overline{E}_{1}, \dotsc ,\overline{E}_{m}^{V} \right)
\]
is an orthonormal frame on $(TM,g^{f})$, where $\overline{E}_{1} = \frac{1}{\sqrt{\alpha}}E_{1}^{V}$ and $\overline{E}_{i} = E_{i}$ for all $i \ge 2$.  From formula (\ref{Eq0.2}), Theorem~\ref{Theo1.1}, Corollary~\ref{Cor.1-2} and Corollary~\ref{Cor_2.1}, we obtain:
\begin{align*}
  \tau(Id_{\alpha}) 
  &= \sum_{i = 1}^{m}\widehat{\nabla}_{E^{H}_{i}}E^{H}_{i} - \sum_{i = 1}^{m}\nabla^{f}_{E^{H}_{i}}E^{H}_{i}+ \sum_{i = 1}^{m}\widehat{\nabla}_{\overline{E}^{V}_{i}}\overline{E}^{V}_{i}-\sum_{i = 1}^{m}\nabla^{f}_{\overline{E}^{V}_{i}}\overline{E}^{V}_{i}\\
  &= -\sum_{i = 1}^{m}\nabla^{f}_{\overline{E}^{V}_{i}}\overline{E}^{V}_{i}\\
  &= -\sum_{i = 1}^{m}\overline{E}_{i}(f)\big{(\nabla_{\overline{E}_{i}}grad(f)\big)}^{H}\\
  &= -\Big[\frac{grad(\alpha)}{2\alpha}\Big]^{H}\\
  &= -\frac{\|grad(f)\|}{1+\|grad(f)\|^{2}}\big{(grad(\|grad(f)\|)\big)}^{H}.
\qedhere
\end{align*}
\end{proof}

\begin{theorem}\label{Th4.4}
Let $(M,g)$ be a flat Riemannian manifold, $f:M \rightarrow ]0,+\infty[$ be a smooth function on $M$ and $(TM,g^{f})$) its tangent bundle equipped with the Mus-Gradient metric. Then the identity $Id_{\alpha}:(TM,g^{f})\mapsto (TM,\widehat{g})$ is biharmonic if and only if
\begin{equation}\label{Eq4.2.1}
    \tau_{2}(Id_{\alpha}) = \sum_{i = 1}^{m}\Big[\nabla_{E_{i}}\nabla_{E_{i}}V -\nabla_{\nabla_{E_{i}}E_{i}}V\big]^{H} = \big[J_{Id_{M}}(W)\big]^{H} = 0,
\end{equation}
where $Id_{M}$ denotes the identity map of the manifold $(M,g)$ and $W = \frac{grad(\alpha)}{\alpha}$.
\end{theorem}

\begin{proof}
Respectively to the notations above, from formula (\ref{Eq0.1}), Theorem~\ref{Theo1.1} and Corollary~\ref{Cor.1-2}, we have:
\begin{align*}
\tau_{2}(Id_{\alpha}) 
&= \sum_{i = 1}^{m}\widehat{\nabla}_{E^{H}_{i}}\widehat{\nabla}_{E^{H}_{i}}W^{H} -\widehat{\nabla}_{\widehat{\nabla}_{E^{H}_{i}}E^{H}_{i}}W^{H}+\widehat{\nabla}_{\overline{E}^{V}_{i}}\widehat{\nabla}_{\overline{E}^{V}_{i}}W^{H} -\widehat{\nabla}_{\widehat{\nabla}_{\overline{E}^{V}_{i}}\overline{E}^{V}_{i}}V^{H}\\
&= \sum_{i = 1}^{m}\widehat{\nabla}_{E^{H}_{i}}\widehat{\nabla}_{E^{H}_{i}}W^{H} -\widehat{\nabla}_{\widehat{\nabla}_{E^{H}_{i}}E^{H}_{i}}W^{H}\\
&= \sum_{i = 1}^{m}\Big[\nabla_{E_{i}}\nabla_{E_{i}}W -\nabla_{\nabla_{E_{i}}E_{i}}W\Big]^{H}\\
&= \big[J_{Id_{\alpha}}(W)\big]^{H}.
\qedhere
\end{align*}
\end{proof}

\begin{lemma}\label{Lem4.2.2}
Let $M = \mathbb{R}^{m}$ be the real euclidean manifold, $f(x_{1},\dots,x_{m}) = f(x_{1})$ be a smooth function strictly positive on $M$ and $(TM,g^{f})$) its tangent bundle equipped with the Mus-Gradient metric. Then the identity $Id_{\alpha}:(TM,g^{f})\mapsto (TM,\widehat{g})$ is biharmonic if and only if
\begin{equation}\label{Eq4.2.3}
    \left( \frac{\alpha'}{2\alpha} \right)'' = 0.
\end{equation}
\end{lemma}

\begin{proof}
Respectively to the notations above, we have:
\begin{align*}
grad(f) &= f',\\
E_{1} &= \frac{grad(f)}{\|grad(f)\|}\partial_{1} = \partial_{1},\\
\overline{E}_{i} &= \partial_{i}\quad,\quad 2\leq i\leq m\\
W &= \frac{\alpha'}{2\alpha}\partial_{1}.
\end{align*}
From Theorem~\ref{Th4.2} we obtain:
\begin{align*}
J_{Id_{\alpha}}(W) 
&= \sum_{i = 1}^{m}\Big[\nabla_{E_{i}}\nabla_{E_{i}}W -\nabla_{\nabla_{E_{i}}E_{i}}W\Big]\\
&= \nabla_{E_{1}}\nabla_{E_{1}}W\\
&= \big(\frac{\alpha'}{2\alpha}\big)''\partial_{1}.
\qedhere
\end{align*}
\end{proof}

According to Lemma~\ref{Lem4.2.2}, we have the following example.

\begin{example}\label{Ex4.2}
If we set
\[
f(x,x_{2},\dots,x_{m}) = \int\sqrt{\exp(ax^{2}+bx+c) -1},
\]
then $f$ is solution of the equation (\ref{Eq4.2.3}) and $Id_{\alpha}$ is proper biharmonic.
\end{example}

\begin{theorem}\label{Th5.1}
Let $(M,g)$ be a flat Riemannian manifold, $f:M \rightarrow ]0,+\infty[$ be a smooth function on $M$ and $(TM,\widetilde{g})$) its tangent bundle equipped with the Mus-Sasaki metric. Then the tension (resp bitension) field of the identity $Id_{\alpha}:(TM,\widetilde{g})\mapsto (TM,\widehat{g})$ is given by
\begin{gather}
    \tau(Id_{f}) = \frac{m}{2f}\big{(grad(f)\big)}^{H}\label{Eq5.2.1},\\
    \tau_{2}(Id_{f}) = \frac{m} {f\sqrt{f}}\big{(\nabla_{grad(f)}grad(\sqrt{f})\big)}^{V}+\frac{m^{2}}{2f^{2}}\big{(\nabla_{grad(f)}grad(f)\big)}^{H}.\label{Eq5.2.2}
\end{gather}
\end{theorem}

\begin{proof}
Let $(E_{1},\dotsc,E_{m})$ be an orthonormal frame on $(M,g)$ such that $\nabla_{E_{i}}E_{j} = 0$ for all $i,j \in \{1, \dotsc, m\}$.  Then $\left( E_{1}^{H},\dotsc,E_{m}^{H},\frac{1}{\sqrt{f}}E_{1}^{V},\dotsc,\frac{1}{\sqrt{f}}E_{m}^{V} \right)$ is an orthonormal frame on $(TM,\widetilde{g})$. From Formula~\ref{Eq0.2}, Corollary~\ref{Cor.0.1} and Lemma~\ref{Lem.0.2}, we obtain:
\begin{align*}
\tau(Id_{f}) 
=& Tr_{\widetilde{g}}(Id_{f}))\\
=& \sum_{i}\Big[\widehat{\nabla}_{E_{i}^{H}}E_{i}^{H} - \widetilde{\nabla}_{E_{i}^{H}}E_{i}^{H}+\widehat{\nabla}_{\frac{1}{\sqrt{f}}E_{i}^{V}}\frac{1}{\sqrt{f}}E_{i}^{V} - \widetilde{\nabla}_{\frac{1}{\sqrt{f}}E_{i}^{V}}\frac{1}{\sqrt{f}}E_{i}^{V}\Big]\\
=& -\sum_{i} \widetilde{\nabla}_{\frac{1}{\sqrt{f}}E_{i}^{V}}\frac{1}{\sqrt{f}}E_{i}^{V}\\
=& \frac{m}{2f}\big{(grad(f)\big)}^{H}, \\
\tau_{2}(Id_{f})
=& \sum_{i}\Big[\widehat{\nabla}_{E_{i}^{H}}\widehat{\nabla}_{E_{i}^{H}}\frac{m}{2f}{(grad(f))}^{H}-
\widehat{\nabla}_{\widetilde{\nabla}_{E_{i}^{H}}E_{i}^{H}}\frac{m}{2f}{(grad(f))}^{H}\Big]\\
&+\sum_{i}\Big[\widehat{\nabla}_{\frac{1}{\sqrt{f}}E_{i}^{V}}\widehat{\nabla}_{\frac{1}{\sqrt{f}}E_{i}^{V}}\frac{m}{2f}{(grad(f))}^{H}-
\widehat{\nabla}_{\widetilde{\nabla}_{\frac{1}{\sqrt{f}}E_{i}^{V}}\frac{1}{\sqrt{f}}E_{i}^{V}}\frac{m}{2f}{(grad(f))}^{H}\Big]\\
=& \sum_{i}\Big[\widehat{\nabla}_{\frac{1}{\sqrt{f}}E_{i}^{V}}\widehat{\nabla}_{\frac{1}{\sqrt{f}}E_{i}^{V}}\frac{m}{2f}{(grad(f))}^{H}-
\widehat{\nabla}_{\widetilde{\nabla}_{\frac{1}{\sqrt{f}}E_{i}^{V}}\frac{1}{\sqrt{f}}E_{i}^{V}}\frac{m}{2f}{(grad(f))}^{H}\Big]\\
=& \sum_{i}\Big[-\widehat{\nabla}_{\widetilde{\nabla}_{\frac{1}{\sqrt{f}}E_{i}^{V}}}\frac{1}{\sqrt{f}}E_{i}^{V}\frac{m}{2f}{(grad(f))}^{H}\Big]\\
=& \widehat{\nabla}_{\frac{m}{2f}{(grad(f))}^{H}}\frac{m}{2f}{(grad(f))}^{H}\\
=& \frac{m^{2}}{2f^{2}}\Big{(\nabla_{grad(f)}grad(f)\Big)}^{H}\\
=& \frac{m^{2}}{4f^{2}}\Big{(grad(f)(\|grad(f)\|^{2})\Big)}^{H}.
\qedhere
\end{align*}
\end{proof}

From Theorem~\ref{Th5.1} we have the following theorem.

\begin{theorem}\label{Th5.2}
Let $(M,g)$ be a flat Riemannian manifold, $f:M \rightarrow ]0,+\infty[$ be a smooth function on $M$ and $(TM,\widetilde{g})$) its tangent bundle equipped with the Mus-Sasaki metric. Then the identity $Id_{\alpha}:(TM,\widetilde{g})\mapsto (TM,\widehat{g})$ is proper biharmonic if and only
\begin{equation}\label{Eq5.1}
    \|grad(f)\| = const>0.
\end{equation}
\end{theorem}

\begin{example}\label{Ex5.1}
Let $M = ]0,+\infty[^{p}\times\mathbb{R}^{m-p}$. If we set
\[
f(x_{1},x_{2},\dots,x_{m}) = \sum_{i = 1}^{p}a^{i}x_{i}+b^{i},\quad a^{i},b^{i}>0.
\]
then $f$ is solution of the equation (\ref{Eq5.1}) and $Id_{f}$ is proper biharmonic.
\end{example}

\begin{theorem}\label{Th5.3}
Let $(M,g)$ be a flat Riemannian manifold, $f:M \rightarrow ]0,+\infty[$ be a smooth function on $M$ and $(TM,\widetilde{g})$) its tangent bundle equipped with the Mus-Sasaki metric. Then the tension (resp bitension) field of the identity $\widehat{Id}_{f}:(TM,\widetilde{g})\mapsto (TM,\widehat{g})$ is given by
\begin{equation}\label{Eq5.2.3}
    \tau(\widehat{Id}_{f}) = -\frac{m}{2}\big{(grad(f)\big)}^{H},
\end{equation}
resp
\begin{equation}\label{Eq5.2.4}
    \tau_{2}(\widehat{Id}_{f}) = -\frac{m^{2}}{8}\big{(grad\|grad(f)\|^{2}\big)}^{H}+\frac{m}{2}\big{(Tr_{g}\nabla^{2}(grad(f))\big)}^{H}.
\end{equation}
\end{theorem}

\begin{proof}
Let $(E_{1},\dotsc,E_{m})$ be an orthonormal frame on $(M,g)$ such that $\nabla_{E_{i}}E_{j} = 0$ for all $i,j \in \{1, \dotsc, m\}$. We have:
\begin{align*}
\tau(\widehat{Id}_{f})
=& Tr_{\widetilde{g}}(Id_{f}))\\
=& \sum_{i}\Big[\widetilde{\nabla}_{E_{i}^{H}}E_{i}^{H}-\widehat{\nabla}_{E_{i}^{H}}E_{i}^{H} +\widetilde{\nabla}_{E_{i}^{V}}E_{i}^{V}-\widehat{\nabla}_{E_{i}^{V}}E_{i}^{V} - \Big]\\
=& -\sum_{i} \widetilde{\nabla}_{E_{i}^{V}}E_{i}^{V}\\
=& -\frac{m}{2}\big{(grad(f)\big)}^{H}, \\
\tau_{2}(\widehat{Id}_{f})
=& \sum_{i}\frac{m}{2}\widetilde{R}(grad{(f)}^{H},E_{i}^{H})E_{i}^{H}+\sum_{i}\frac{m}{2}\widetilde{R}(grad{(f)}^{H},E_{i}^{V})E_{i}^{V}\\
&+\sum_{i}\Big[\widetilde{\nabla}_{E_{i}^{H}}\widetilde{\nabla}_{E_{i}^{H}}\frac{m}{2}{(grad(f))}^{H}-
\widetilde{\nabla}_{\widehat{\nabla}_{E_{i}^{H}}E_{i}^{H}}\frac{m}{2}{(grad(f))}^{H}\Big]\\
&+\sum_{i}\Big[\widetilde{\nabla}_{E_{i}^{V}}\widetilde{\nabla}_{E_{i}^{V}}\frac{m}{2}{(grad(f))}^{H}-
\widetilde{\nabla}_{\widehat{\nabla}_{E_{i}^{V}}E_{i}^{V}}\frac{m}{2}{(grad(f))}^{H}\Big]\\
=& \sum_{i}\frac{m}{2}\widetilde{R}(grad{(f)}^{H},E_{i}^{V})E_{i}^{V}
+\frac{m}{2}\sum_{i}\widetilde{\nabla}_{E_{i}^{H}}\widetilde{\nabla}_{E_{i}^{H}}{(grad(f))}^{H}\\
&+\frac{m}{2}\sum_{i}\widetilde{\nabla}_{E_{i}^{V}}\widetilde{\nabla}_{E_{i}^{V}}{(grad(f))}^{H}\\
=& -\big{(\frac{m^{2}}{4}\nabla_{grad(f)}grad(f)\big)}^{H}+\frac{m^{2}}{8f}\|grad(f)\|^{2}\big{(grad(f)\big)}^{H}\\
&+\frac{m}{2}\sum_{i}\big{(\nabla_{E_{i}}\nabla_{E_{i}}(grad(f))\big)}^{H}+\frac{m}{4}\sum_{i}\widetilde{\nabla}_{E_{i}^{V}}(grad(f)(f))E_{i}^{V}\\
=& -\frac{m^{2}}{8}\big{(grad\|grad(f)\|^{2}\big)}^{H}+\frac{m^{2}}{8f}\|grad(f)\|^{2}\big{(grad(f)\big)}^{H}\\
&+\frac{m}{2}\sum_{i}\big{(\nabla_{E_{i}}\nabla_{E_{i}}(grad(f))\big)}^{H}-\frac{m^{2}}{8}\|grad(f)\|^{2}{(grad(f))}^{H}\\
=& -\frac{m^{2}}{8}\big{(grad\|grad(f)\|^{2}\big)}^{H}+\frac{m}{2}\big{(Tr_{g}\nabla^{2}(grad(f))\big)}^{H}.
\qedhere
\end{align*}
\end{proof}

By Theorem~\ref{Th5.3} we obtain the following theorem.

\begin{theorem}\label{Th5.4}
Let $(M,g)$ be a flat Riemannian manifold, $f:M \rightarrow ]0,+\infty[$ be a smooth function on $M$ and $(TM,\widetilde{g})$) its tangent bundle equipped with the Mus-Sasaki metric. Then the identity $\widehat{Id}_{f}:(TM,\widehat{g})\mapsto (TM,\widetilde{g})$ is a proper biharmonic if and only
\[
\frac{m}{4}grad(\|grad(f)\|^{2}) = Tr_{g}\nabla^{2}(grad(f))\;.
\]
\end{theorem}

\begin{lemma}\label{Lem5.1.1}
Let $M = ]0,+\infty[\times\mathbb{R}^{m-1}$ be the real euclidean manifold, $f: M \to \mathbb R$ be the function $f(x_{1},\dotsc,x_{m}) = f(x_{1})$, which is strictly positive, and $(TM,\widetilde{g})$ be its tangent bundle equipped with the Mus-Sasaki metric. Then the identity $\widehat{Id}_{f}:(TM,\widehat{g})\mapsto (M,\widetilde{g})$ is biharmonic if and only if
\begin{equation}\label{Eq5.2.5}
    f''' = \frac{m}{4}\big({(f')}^{2}\big)'.
\end{equation}
\end{lemma}

\begin{example}[proper biharmonic map]\label{Ex5.2}
Let $M = ]0,+\infty[^{p}\times\mathbb{R}^{m-p}$. If we set
\[
f(x_{1},x_{2},\dots,x_{m}) = \frac{4}{m}\ln(mx+c),
\]
then $f$ is solution of the equation (\ref{Eq5.2.5}) and $\widehat{Id}_{f}$ is proper biharmonic.
\end{example}

\subsection*{Acknowledgments}

This note was supported by L.G.A.C.A. Laboratory of Saida university and P.R.F.U project. The authors would like to thank Professor Mustapha Djaa for his assistance in the idea of this work and its guidelines, as well the anonymous referee.

\EditInfo{April 22, 2021}{December 21, 2021}{Haizhong Li}

\end{paper}
\begin{references}

\refer{Paper}{Aba}
\Rauthor{Abbassi M.T.K. and Sarih M.}
\Rtitle{On natural metrics on tangent bundles of Riemannian manifolds}
\Rjournal{Arch. Math}
\Rvolume{41}
\Ryear{2005}
\Rpages{71-92}

\refer{Paper}{BFS}
\Rauthor{Baird P., Fardoun A. and Ouakkas S.}
\Rtitle{Conformal and semi-conformal biharmonic maps}
\Rjournal{Ann. Global Anal. Geom}
\Rvolume{34}
\Ryear{2008}
\Rpages{403-414}

\refer{Paper}{BK}
\Rauthor{P. Baird and D. Kamissoko}
\Rtitle{On constructing biharmonic maps and metrics}
\Rjournal{Ann. Global Anal. Geom}
\Rvolume{23}
\Ryear{2003}
\Rpages{65-75}

\refer{Paper}{AB}
\Rauthor{Balmus A.}
\Rtitle{Biharmonic properties and conformal changes}
\Rjournal{An. Stiint. Univ. Al.I.Cuza Iasi Mat. (N.S)}
\Rvolume{50}
\Ryear{2004}
\Rpages{367-372}

\refer{Paper}{BMO}
\Rauthor{A. Balmus, S. Montaldo and C. Oniciuc}
\Rtitle{Biharmonic maps between warped product manifolds}
\Rjournal{J. Geom. Phys}
\Rvolume{57}
\Ryear{2008}
\Rpages{449-466}

\refer{Paper}{C.M.O}
\Rauthor{Caddeo R., Montaldo S. and Oniciuc C.}
\Rtitle{Biharmonic submanifolds of $\mathbb{S}^{3}$}
\Rjournal{Int. J. Math}
\Rvolume{12}
\Ryear{2001}
\Rpages{867-876}

\refer{Paper}{Key5}
\Rauthor{Cheeger J. and Gromoll D.}
\Rtitle{On the structure of complete manifolds of nonnegative curvature}
\Rjournal{Ann. of Math}
\Rvolume{2}
\Ryear{1972}
\Rnumber{96}
\Rpages{413-443}

\refer{Proceedings}{D.L}
\Rauthor{Djaa M. and Latti F.}
\Rtitle{Non linear analysis on manifolds On generalized $F$-energy variation between Riemannian manifolds}
\Reditor{AIP Conference Proceedings}
\Rjournal{2074}
\Rpublisher{2019}
\Rpages{1-13}

\refer{Paper}{Key3}
\Rauthor{Dombrowski P.}
\Rtitle{On the structure of complete manifolds of nonnegative curvature}
\Rjournal{Journal für die reine und angewandte Mathematik}
\Rvolume{210}
\Ryear{1962}
\Rpages{73-88}

\refer{Paper}{J.L1}
\Rauthor{Eells J.and Lemaire L.}
\Rtitle{A report on harmonic maps}
\Rjournal{Bull. London Math. Soc}
\Rvolume{10}
\Ryear{1978}
\Rpages{1-68}

\refer{Paper}{J.L2}
\Rauthor{Eells J.and Lemaire L.}
\Rtitle{Another report on harmonic maps}
\Rjournal{Bull. London Math. Soc}
\Rvolume{20}
\Ryear{1988}
\Rpages{385-524}

\refer{Paper}{Key19}
\Rauthor{Gezer A.}
\Rtitle{On the tangent bundle with deformed Sasaki metric}
\Rjournal{Int. Electron. J. Geom}
\Rvolume{6}
\Ryear{2013}
\Rnumber{2}
\Rpages{19-31}

\refer{Paper}{Key9}
\Rauthor{Gudmunsson S. and Kappos E.}
\Rtitle{On the Geometry of Tangent Bundles}
\Rjournal{Expo. Math}
\Rvolume{20}
\Ryear{2002}
\Rpages{1-14}

\refer{Paper}{L.D.Z}
\Rauthor{Latti F., Djaa M. and Zagane A.}
\Rtitle{Mus-Sasaki Metric and Harmonicity}
\Rjournal{Mathematical Sciences and Applications E-Notes}
\Rvolume{6}
\Ryear{2018}
\Rnumber{1}
\Rpages{29-36}

\refer{Paper}{Lu}
\Rauthor{Lu W. J.}
\Rtitle{Geometry of warped product manifolds and its five applications}
\Rjournal{Ph. D. thesis, Zhejiang University}
\Ryear{2013}

\refer{Paper}{MC.D}
\Rauthor{Mohammed Cherif A. and Djaa M.}
\Rtitle{Harmonic Maps and Torse-Forming Vector Fields}
\Rjournal{International Electronic Journal of Geometry}
\Rvolume{13}
\Ryear{2020}
\Rnumber{1}
\Rpages{87-93}

\refer{Paper}{Key8}
\Rauthor{Musso E. and Tricerri F.}
\Rtitle{Riemannian metrics on tangent bundles}
\Rjournal{Ann.Mat.Purz Appl}
\Rvolume{4}
\Ryear{1988}
\Rnumber{150}
\Rpages{1-10}

\refer{Paper}{O.N.D}
\Rauthor{Ouakkas S., Nasri R. and Djaa M.}
\Rtitle{On the f-harmonic and f-biharmonic maps}
\Rjournal{J. P. Journal of Geom. and Top}
\Rvolume{10}
\Ryear{2010}
\Rnumber{1}
\Rpages{11-27}

\refer{Paper}{KZD}
\Rauthor{Kada Benothmane A., Zagane A. and Djaa M.}
\Rtitle{On Generalized Cheeger-Gromoll Metric And Harmonicity}
\Rjournal{Commun. Fac. Sci. Univ. Ank. Ser. A1 Math. Stat}
\Rvolume{69}
\Ryear{2020}
\Rnumber{1}
\Rpages{629-645}

\refer{Paper}{Key15}
\Rauthor{Salimov A. A., Gezer A. and Akbulut K.}
\Rtitle{Geodesics of Sasakian metrics on tensor bundles}
\Rjournal{Mediterr. J. Math}
\Rvolume{6}
\Ryear{2009}
\Rnumber{2}
\Rpages{135-147}

\refer{Paper}{Key16}
\Rauthor{Salimov A. A. and Gezer A.}
\Rtitle{On the geometry of the (1, 1)-tensor bundle with Sasaki type metric}
\Rjournal{Chinese Annals of Mathematics, Series B}
\Rvolume{32}
\Ryear{2011}
\Rnumber{3}
\Rpages{369-386}

\refer{Paper}{Key17}
\Rauthor{Salimov A. A. and Agca F.}
\Rtitle{Some Properties of Sasakian Metrics in Cotangent Bundles}
\Rjournal{Mediterranean Journal of Mathematics}
\Rvolume{8}
\Ryear{2011}
\Rnumber{2}
\Rpages{243-255}

\refer{Paper}{Key18}
\Rauthor{Salimov A. A. and Kazimova S.}
\Rtitle{Geodesics of the Cheeger-Gromoll Metric}
\Rjournal{Turk J Math}
\Rvolume{33}
\Ryear{2009}
\Rpages{99-105}

\refer{Paper}{Key2}
\Rauthor{Sasaki S.}
\Rtitle{On the differential geometry of tangent bundles of Riemannian manifolds}
\Rjournal{Tohoku Math. J}
\Rvolume{10}
\Ryear{1958}
\Rpages{338-354}

\refer{Arxiv}{Key4}
\Rauthor{Wang J. and Yong Wang}
\Rtitle{On the Geometry of Tangent Bundles with the Rescaled Metric}
\Rarxivid{arXiv:1104.5584v1 [math.DG] 29 Apr 2011.}

\refer{Book}{Key7}
\Rauthor{Yano K. and Ishihara S.}
\Rtitle{Tangent and Cotangent Bundles}
\Rpublisher{Marcel Dekker. INC. New York}
\Ryear{1973}

\refer{Paper}{Z.D}
\Rauthor{Zagane A. and Djaa M.}
\Rtitle{Geometry of Mus-Sasaki metric}
\Rjournal{Communication in Mathematics}
\Rvolume{26}
\Ryear{2018}
\Rpages{113-126}

\end{references}
